\documentclass[12pt]{article}
\usepackage{amsmath,amssymb}
\usepackage{latexsym}
\usepackage{epsfig}

\usepackage{psfrag}

\begin{document}
\title
 { ~ \\ ~ \\ Unitary matrices associated with Butson-Hadamard matrices}

\author{Alberto Mart\'{i}n M\'endez$^1$}

\maketitle

\footnotetext[1] {Department of Applied Mathematics II, University
of Vigo, Vigo, Spain \\
email:amartin@dma.uvigo.es}
\begin{abstract}
\noindent We give an example of a $BH(5,5)$ matrix $M$ satisfying that the eigenvalues of the unitary matrix $\frac{1}{\sqrt{5}}M$ are all primitive $10$-th roots of unity and such that $\frac{1}{5}M^3$ is not a $BH(5,5)$ matrix. This example gives a negative answer to a conjecture proposed by R. Egan and P. \'O Cath\'ain.

\end{abstract}

\vspace{1cm}
\noindent A Butson-Hadamard matrix is an $m\times m$ complex matrix $M$ whose entries are all $l$-th roots of unity and satisfying $MM^*=mI$, where $M^*$ denotes the conjugate transpose of $M$. These matrices are denoted by $BH(m,l)$~\cite{butson}.

If $M$ is a $BH(m,l)$ matrix, then the matrix $B=\frac{1}{\sqrt{m}}M$ is an unitary matrix. Reciprocally, if $B$ is an $m\times m$ unitary matrix and all the entries of $M=\sqrt{m}B$ are $l$-th roots of unity, then $M$ is a $BH(m,l)$ matrix. We will say that $B=\frac{1}{\sqrt{m}}M$ is the unitary matrix associated with $M$.

In~\cite{conjecture}, pag. 84, R. Egan and P. \'O Cath\'ain pose, by way of conjecture (in their own words, "we feel it should have a positive answer"), the following question, which we first reproduce as stated therein:

\vspace{0.5cm}

{\em If the eigenvalues of $M\in BH(m,l)$ are all primitive $k$-th roots of unity, is it true that $\sqrt{m}^{1-i}M^{i} \in BH(m,l)$ for all $i$ coprime to $k$?}

\vspace{0.5cm}

As it stands, the question contains an inaccuracy, since the condition $MM^*=mI_m$ implies that the modulus of all the eigenvalues of $M$ is $\sqrt{m}$; therefore, none of them can be a (primitive) $k$-th root of unity. Relying on the rest of the article and on the first example below, we can correctly restate the above question as:

\vspace{0.5cm}

{\em If $M\in BH(m,l)$ and the eigenvalues of the associated unitary matrix $B=\frac{1}{\sqrt{m}}M$ are all primitive $k$-th roots of unity, is it true that $\sqrt{m}^{1-i}M^{i} \in BH(m,l)$ for all $i$ coprime to $k$?}

\vspace{0.5cm}

On the order hand, if $M\in BH(m,l)$ and $i\geq 1$ we have that
$$(\sqrt{m}^{1-i}M^{i})(\sqrt{m}^{1-i}M^{i})^*=\sqrt{m}^{2-2i}M^{i}M^{*i}=m^{1-i}m^{i}I_m=mI_m.$$

Therefore, the matrix $\sqrt{m}^{1-i}M^{i}$ belongs to $BH(m,l)$ if and only if all its entries are $l$-th roots of unity. This allows us to reformulate the former question as follows:

\vspace{0.5cm}

{\em If $M\in BH(m,l)$ and the eigenvalues of the associated unitary matrix $B=\frac{1}{\sqrt{m}}M$ are all primitive $k$-th roots of unity, is it true that for all $i$ coprime to $k$ the entries of $\sqrt{m}^{1-i}M^{i}$ are all $l$-th roots of unity?}

\vspace{0.5cm}

Indeed, in the example
$$M=\left( \begin{array} {rr}
1 & 1 \\
i & -i \end{array} \right) \in BH(2,4),$$

\noindent taken from~\cite{conjecture}, the associated unitary matrix $B=\frac{1}{\sqrt{2}}M$ has as eigenvalues $e^{\frac{\pi }{12}i}$ and $e^{\frac{17\pi }{12}i}$; both are primitive $24$-th roots of unity. We are in the hypotheses of our statement for $m=2$, $l=4$ and $k=24$. The first integer coprime to $k$ is $i=5$; the matrix
$$\sqrt{2}^{1-5}M^5=\frac{1}{4}M^5=\left( \begin{array}{rr} i & 1 \\
i & -1 \end{array} \right)$$

\noindent satisfies that all its entries are $4$-th roots of unity, and 
$\frac{1}{4}M^5\in BH(2,4)$.

However, the answer to the conjecture presented by Egan and \'O Cath\'ain is negative. Let us consider the matrix
$$M=\left( \begin{array}{ccccc}
e^{\frac{2\pi}{5}i} & e^{\frac{6\pi}{5}i} & e^{\frac{8\pi}{5}i} & e^{\frac{8\pi}{5}i} & e^{\frac{6\pi}{5}i} \\
e^{\frac{6\pi}{5}i} & e^{\frac{2\pi}{5}i} & e^{\frac{6\pi}{5}i} & e^{\frac{8\pi}{5}i} & e^{\frac{8\pi}{5}i} \\
e^{\frac{8\pi}{5}i} & e^{\frac{6\pi}{5}i} & e^{\frac{2\pi}{5}i} & e^{\frac{6\pi}{5}i} & e^{\frac{8\pi}{5}i} \\
e^{\frac{8\pi}{5}i} & e^{\frac{8\pi}{5}i} & e^{\frac{6\pi}{5}i} & e^{\frac{2\pi}{5}i} & e^{\frac{6\pi}{5}i} \\
e^{\frac{6\pi}{5}i} & e^{\frac{8\pi}{5}i} & e^{\frac{8\pi}{5}i} & e^{\frac{6\pi}{5}i} & e^{\frac{2\pi}{5}i}  \end{array} \right),$$

\noindent which is a circulant symmetric matrix, whose entries are all (primitive) $5$-th roots of unity, and satisfies $MM^*=5I_5$; i.e., $M\in BH(5,5)$. Let us note that $M$ is also an unreal matrix in the sense of~\cite{unreal}. Since the associated unitary matrix $B=\frac{1}{\sqrt{5}}M$ is also a circulant matrix, its eigenvalues (see~\cite{eigenvalues}) are given by $h(\xi ^{j})$, $j=0,1,2,3,4$, where $\xi =e^{\frac{2\pi}{5}i}$ and 
$$h(x)= e^{\frac{2\pi}{5}i}+ e^{\frac{6\pi}{5}i}x+ e^{\frac{8\pi}{5}i}x^2+ e^{\frac{8\pi}{5}i}x^3+ e^{\frac{6\pi}{5}i}x^4.$$

A  calculation shows that the eigenvalues of $B$ are $e^{\frac{3\pi}{5}i}$, $e^{\frac{3\pi}{5}i}$, $e^{\frac{\pi}{5}i}$, $e^{\frac{\pi}{5}i}$ and $e^{\frac{7\pi}{5}i}$. All of them are primitive $10$-th roots of unity. So, we are again in the hypotheses of the statement for $m=l=5$ and $k=10$. For $i=3$, which is coprime to $k$, we have that the matrix
$$\sqrt{5}^{1-3}M^3=\frac{1}{5}M^3$$

\noindent has only three different entries, $e^{\frac{3\pi}{5}i}$, $e^{\frac{\pi}{5}i}$ (both of which coincide with eigenvalues of $B$) and $e^{\frac{9\pi}{5}i}$. None of these entries is a $5$-root of unity; thus, the answer to the raised question is negative.

Finally, we can observe that in the two examples above the entries of the computed matrices $\sqrt{m}^{1-i}M^{i}$ are all $k$-th roots of unity. Nevertheless, it should be noted that, in general, this is also untrue. For example, take the Hadamard matrix with constant diagonal
$$M=\left( \begin{array}{rrrr}
-1 & 1 & -1 & 1 \\
-1 & -1 & -1 & -1 \\
1 & 1 & -1 & -1  \\
-1 & 1 & 1 & -1 \end{array} \right) \in BH(4,2).$$

The minimal polynomial of $M$ is $p(x)=x^2+2x+4$ (in fact, we have $(H+I)^t=-H-I$, i.e., $H=-2I-H^t$, and plugging it into $p(x)$ we obtain $p(H)=0$) and then, the eigenvalues of the associated unitary matrix $B=\frac{1}{2}M$ are $e^{\frac{2\pi}{3}i}$, $e^{\frac{2\pi}{3}i}$, $e^{\frac{4\pi}{3}i}$ and $e^{\frac{4\pi}{3}i}$, all of them primitive cube roots of unity. Taking $i=2$, which is coprime to $k=3$, we obtain
$$\sqrt{4}^{1-2}M^2=\frac{1}{2}M^2=\left( \begin{array}{rrrr}
-1 & -1 & 1 & -1 \\
1 & -1 & 1 & 1 \\
-1 & -1 & -1 & 1 \\
1 & -1 & -1 & -1 \end{array} \right),$$

\noindent and the entries of $\frac{1}{2}M^2$ are all square roots of unity, not cube roots of unity.

\end{document}